
\input amstex
\TagsOnRight
\raggedbottom
\hoffset=0cm\voffset=0cm
\magnification=\magstep1

\baselineskip=15pt

\def\ref#1{\noindent\item{[#1.]}}
\def\spa#1{\noindent\item {#1}}

\vskip 9cm
\hfill {\tenrm NUS/HEP/95201}
\par\noindent
\hfill {\tenrm q-alg/9508015}
\vskip 2cm
\centerline {\bf Real Forms of the Oscillator Quantum Algebra}
\centerline {\bf and its Representations}
\vskip.7cm
\centerline {C.H. Oh $^{\dag}$ and K. Singh }
\vskip .5cm
{\sl {\centerline {Department of Physics}}}
{\sl {\centerline {Faculty of Science}}}
{\sl {\centerline {National University of Singapore}}}
{\sl {\centerline {Lower Kent Ridge, Singapore 0511}}}
{\sl {\centerline {Republic of Singapore}}}
\vskip 2cm
{\bf {\centerline {Abstract}}}
\vskip .5cm
\noindent
We consider the conditions under which the $q$-oscillator algebra
becomes a Hopf $*$-algebra. In particular, we show that there are
at least two real forms associated with the algebra. Furthermore, through
the representations, it is shown that they are related to
${\text {su}}_{q^{1/2}}(2)$ with different conjugations.
\vskip 5truecm
\hrule width6cm
\vskip .3truecm
{\sevenrm $^{\dag}$ E-mail: PHYOHCH\@NUSVM.NUS.SG}
\eject
\rm
\par
The recent focus on a class of algebras called quantum algebras, or
sometimes loosely called quantum groups, can be attributed to the
fact that it appears to be the basic algebraic structure underlying
many of the physical theories of current interest. Indeed, since its
emergence from the study of the inverse scattering problem [1], it
is now well established that it is also deeply rooted in other areas such
as exactly soluble models [2], factorizable $S$-matrix theory [3] and
conformal field theory [4].
Mathematically, it has been shown that these structures can best be
decribed by a general class of associative algebras called Hopf
algebras [5].
\par
To date many explicit examples of these algebras have been
furnished. In particular, the quantum algebra ${\text {su}}_q(2)$, by
virtue of its simplicity and relevance to physical systems, has
attracted a lot of attention. In fact it has now, more or less, assumed a
paradigmatic role in the study of quantum algebras. Of no less
importance are the oscillator algebras, or in this case
$q$-oscillators. First introduced by Macfarlane [6] and
independently by Biedernharn [7], they have been used in giving
a realization of ${\text {su}}_q(2)$.
\par
More recently it has been shown [8] (see also [9]) that the $q$-oscillator
algebra, when expressed in a commutator form, may itself support a Hopf
structure. In other words, it itself is a quantum algebra. In the
following we examine this algebra under the conditions of a Hopf
$*$-algebra [10-12]. This essentially provides the {\it real forms}
of the algebra or alternatively allows for the restriction of the
algebra over the real number field since in physics one is often
interested in operators that are hermitian.
Here, we show that for the conjugation as implied in the algebra
of ref.[8], the only admissable values of $q$ are  $\{q\in {\bold
C}; \vert q\vert$=1\}. We also present an alternative
conjugation which allows for real $q$. In both cases, it is shown,
through their representations that they are closely related to
${\text {su}}_{q^{1/2}}(2)$.
\par
We begin by considering the associative algebra ${\Cal H}_q(1)$
generated by the elements
$a$, $\overline a$ and $N$ satisfing [8]
$$\align
&[a,{\overline a}]=[N+1]_q-[N]_q \tag 1a \\
&[N,{\overline a}]={\overline a}\qquad [N,a]=-a \tag 1b\\
\endalign$$
where
$$[x]_q=\frac {q^x-q^{-x}}{q-q^{-1}}. \tag 2$$
It should be noted that we differ from the notation of ref.[8] in two respects.
Firstly, the right hand side of (1a) as given there is $[N+\frac {1}{2}]_q-
[N-\frac {1}{2}]_q$ instead of the one given above. Here we have
made the replacement $N\to N+\frac {1}{2}$. Secondly,
${\overline a}$ instead of $a^+$ has been used, as we will also be
considering the
situation in which $a$ and ${\overline a}$ are not conjugates of
one another. The Hopf structure associated with
${\Cal H}_q(1)$ is given by
$$\align
&\Delta(a) = a\otimes q^{1/2(N+\gamma)} + q^{-1/2(N+\gamma)}\otimes a
\tag 3a \\
&\Delta({\overline a}) = {\overline a}\otimes q^{1/2(N+\gamma )}
+ q^{-1/2(N+\gamma )}\otimes {\overline a} \tag 3b \\
&\Delta(N) = N\otimes {\bold l}
+ {\bold l}\otimes N +\gamma {\bold l}\otimes {\bold l} \tag 3c \\
&\epsilon (a) =\epsilon (\overline a) = 0\quad
\epsilon (N)= -\gamma \tag 3d \\
&S(a)= -q^{-1/2}a\quad S({\overline a})=-q^{1/2}{\overline a}\quad
S(N)=-N-2\gamma {\bold l} \tag 3e \\
\endalign$$
where, with $q=e^{\epsilon}$,
$$\gamma =\frac {1}{2}-i\frac {(2l +1)\pi}{2\epsilon} \qquad
l\in {\bold Z}. \tag 4$$
Here $\Delta$, $\epsilon$ and $S$ are the coproduct,
counit and antipode respectively.
\par
Next, let us consider its conjugations. In the context of Hopf algebras,
this requires the notion of a $*$-algebra which is also compatible with
the underlying Hopf structure. More precisely, one defines a Hopf $*$-algebra
as a Hopf algebra ${\Cal A}$ which is equipped with an involution
$a \mapsto a^+$ of ${\Cal A}$ into itself such that
$$\alignat3
&(\alpha a+\beta b)^+={\overline \alpha} a^+ +\overline {\beta}b^+ &&\qquad
(ab)^+=b^+a^+ &&\qquad (a^+)^+ = a \tag5a \\
&\Delta (a^+)=\Delta (a)^+&&\qquad
\epsilon(a^+)={\overline {\epsilon (a)}}&&\qquad S(S(a)^+)^+=a\tag 5b\\
\endalignat$$
for any $a,b\in {\Cal A}$ and $\alpha ,\beta \in {\bold C}$ [10].
It should be noted that the above involution is sometimes refered to
as the {\it standard conjugation}. It is also possible to introduce
yet another involution, called the {\it non-standard
conjugation}. Here the action on $\Delta$ reads as [11,12]
$$\Delta (a)^+=\Delta '(a^+)\tag 6a$$
where $\Delta '$ is the opposite coproduct while the antipode
satisfies
$$S(a^+)=S(a)^+.\tag 6b$$
\par
Now returning to ${\Cal H}_q(1)$, let us see whether the canonical conjugation
$$ (a)^+={\overline a}\qquad ({\overline a})^+=a\qquad (N)^+=N \tag
7$$
furnished in ref.[8], satisfies the necessary requirements.
First as a $*$-algebra, (5a) requires that the conjugation
be compatible with the commutation relations. Explicitly, by conjugating both
sides of (1a) and using (7), one has
$$[a,{\overline a}]=[N+1]_{\overline q}-[N]_{\overline q}\tag 8$$
where ${\overline q}$ denotes the complex conjugate of $q$. It is
clear then, that for the relations to be preserved we must restrict
$q$ to be either real or $\vert q\vert=1$ for $q\in {\bold C}$. As
for the compatibility with the Hopf structure, one finds that only
the $\vert q\vert =1$ case is admissable.
To see that it fails for $q\in {\bold R}$, it suffices to
consider conditions (5b) or (6a) when applied to $N$. With $N^+=N$, we have
$$\Delta (N^+)=\Delta (N)=N\otimes {\bold l}+{\bold l}\otimes N +
\gamma {\bold l}\otimes {\bold l}\tag 9a$$
while
$$\align
\Delta (N)^+&=(N\otimes {\bold l}+{\bold l}\otimes N +
\gamma {\bold l}\otimes {\bold l})^+ \\
 &=N\otimes {\bold l}+{\bold l}\otimes N +
{\overline \gamma} {\bold l}\otimes {\bold l} \tag 9b
\endalign$$
from which we surmise that $\Delta (N^+)\ne \Delta (N)^+$ since
$\gamma \ne {\overline \gamma}$. It is also
easy to see that the non-standard conjugation does not hold either.
For the $\vert q\vert =1$ case , we note that $\gamma$ (eqn.(4) with
$\epsilon \to i\epsilon$) now reads as
$$\gamma =\frac {1}{2}-\frac {(2l +1)\pi}{2\epsilon} \qquad
l\in {\bold Z} \tag 10$$
and is real. In this case
both the standard as well as the non-standard conjugation holds for $N$.
For $a$, however, only the non-standard conjugation is admissable. This is
evident
from the following:
$$\align
\Delta (a^+)&= a^+\otimes q^{1/2(N+\gamma)} + q^{-1/2(N+\gamma)}\otimes a^+
\tag 11a \\
\Delta '(a^+)&= a^+\otimes q^{-1/2(N+\gamma)} + q^{1/2(N+\gamma)}\otimes a^+
\tag 11b \\
\Delta (a)^+&= a^+\otimes q^{-1/2(N+\gamma)} + q^{1/2(N+\gamma)}\otimes a^+
.\tag 11c \endalign $$
In addition, the conditions on the counit and the antipode are also
satisfied by all the generators. Thus we
are left to conclude ${\Cal H}_q(1)$ under the involution (7) is a
{\it bonafide} Hopf $*$-algebra only for $\vert q\vert =1$. This
naturally raises the question as to whether other conjugations exist
for which $q\in {\bold R}$ is admissable. In other words, are there
other real forms associated with ${\Cal H}_q(1)$.
\par
To this end we generalize (7) and write instead,
$$(a)^+=\alpha {\overline a}\qquad (\overline a)^+=\beta a \qquad
(N)^+=N+\eta\tag 12 $$
where $\alpha , \beta $ and $\eta $ are constants to be determined.
\plainfootnote {$^{\dag}$}{In (12) above, it is implicit that
$\eta$ means $\eta {\bold l}$.}
By assuming that $q$ is real, it easy to see that the Hopf
$*$-algebra requirements allow for the following involutions:
$$(a)^+=\pm i{\overline a}\qquad (\overline a)^+=\pm i a \qquad
(N)^+=N - i\frac {(2l +1)\pi}{\epsilon}\tag 13 $$
where $l\in {\bold Z}$ here is the same integer as that in (4). It
is important to note that the above involution corresponds to the
standard conjugation in contrast to the non-standard conjugation
for (7).
\par
To obtain a deeper insight into these conjugations, it is
instructive to look at the representations of (1) under these conditions.
Let us first consider the algebra with the involution as defined by (7).
To obtain its representations, we first note that
$N$ commutes with both ${\overline a}a$ and $a{\overline a}$. As a
consequence it is possible to construct a vector $\vert \psi_0>$
which is a simultaneous eigenstate of $N$ and ${\overline a}a$ i.e.
$$N\vert \psi_0>=\nu_0\vert \psi_0>\qquad
{\overline a}a\vert \psi_0>=\lambda_0\vert \psi_0>\tag 14 $$
where $\nu_0$ and $\lambda_0$ are the corresponding eigenvalues.
Note that these values are real since both the operators $N$ and
${\overline a}a$ are hermitian.
Now, from $\vert \psi_0>$ one can construct other eigenstates of $N$
by setting
$$\vert \psi_n>={\overline a}^n\vert \psi_0>\qquad
\vert \psi_{-n}>=a^n\vert \psi_0>\tag 15 $$
for a positive integer $n$. It is easy to verify that
$$N\vert \psi_{\pm n}>=(\nu_0\pm n)\vert \psi_{\pm n}>.\tag 16 $$
Moreover ${\overline a}$ and $a$ act as raising and lowering
operators on the states $\vert \psi_n>$. Indeed, to see this,
consider the operator identities
$$\align
a{\overline a}^n &={\overline a}^na + [n]_{q^{\frac {1}{2}}}
\frac{q^{(N-\frac {n}{2}+1)} +
q^{-(N-\frac {n}{2}+1)}}{q^{\frac {1}{2}} +q^{-\frac {1}{2}}}
{\overline a}^{n-1} \tag 17a \\
{\overline a}{a}^n &=a^n{\overline a} - [n]_{q^{\frac {1}{2}}}
\frac{q^{(N+\frac {n}{2})} +
q^{-(N+\frac {n}{2})}}{q^{\frac {1}{2}} +q^{-\frac {1}{2}}}
{a}^{n-1} \tag 17b
\endalign$$
which can shown by induction. Then from (15)
and the above relations, we have for $n>0$,
$$\alignat2
&{\overline a}\vert \psi_n>=\vert \psi_{n+1}>&&\qquad
a\vert \psi_n>=\lambda_n\vert \psi_{n-1}> \tag 18a \\
&{\overline a}\vert \psi_{-n}>=\lambda_{-n+1}\vert \psi_{-n+1}> &&\qquad
a\vert \psi_{-n}>=\vert \psi_{-n-1}>\tag 18b \\
\endalignat$$
where
$$\lambda_n= \lambda_0+[n]_{q^{{1}\over{2}}}{{q^{\nu_0+{{n}\over{2}}}+
q^{-\nu_0-{{n}\over{2}}}}\over{q^{{1}\over{2}}+q^{-{{1}\over{2}}}}}.\tag
18c$$
It is clear from above that the $q$-oscillator algebra admits an infinite
tower of states, which is a standard feature of the usual oscillator
algebra. Let us now examine the feasibility of truncating these
states. Here as in the trivial case ($q=1$), we can
set $a\vert \psi_0>=0$ without any obstruction. This will mean that
$\lambda_0 =0$ and $\vert \psi_{-n}>=0$ for $n>0$.
Next let us see whether, it is possible to set ${\overline a} \vert
\psi_k>=0$ for some state $\vert \psi_k>$ with integer $k\ge 0$. The
existence of such a state will essentially
limit the non-trivial states to a finite number. Proceeding along
this, we first consider the casmir
$${\Cal C_2} = {\overline a}a-[N]_q =
a{\overline a}-[N+1]_q \tag 19$$
associated with algebra (1)
and assume that $\vert \psi_k>$ does exist. Then by applying
both forms of ${\Cal C}_2$ to $\vert \psi_k>$ and remembering that
$\vert q\vert =1$ with $q=e^{i\epsilon}$ we have
$${\sin}(\epsilon(\nu_0+k+1))={\sin}(\epsilon \nu_0)\tag 20a$$
or equivalently to
$${\cos}(\epsilon (2\nu_0+k+1)/2)~{\sin}(\epsilon (k+1)/2)=0.\tag 20b$$
Now the sine term is non zero for $k\ge 0$ for arbitrary $\epsilon $
while the cosine term can vanish if
$$\epsilon (2\nu_0+k+1) = (2l+1)\pi\qquad l\in {\bold Z}\tag 21$$
holds. Thus if $\nu_0$, which at this stage is still arbitrary, is
suitably chosen such that the above condition is satisfied for some $k$,
then it is possible to have a finite representation.
Implementing this, we have from (21)
$$\nu_0=\frac {-k-1}{2}+\frac{(2l+1)\pi}{2\epsilon}.\tag 22$$
As a consequence, we have
\plainfootnote {$^{\dag}$}
{Here we have included the index $k$ in the state to reflect the
dependence of the representation on this parameter.}
$${\overline a}\vert \psi_{k,n}>=\vert \psi_{k,n+1}>\qquad
a\vert \psi_{k,n}>=\lambda_n\vert \psi_{k,n-1}>\tag 23a $$
with the finiteness in the number of states controlled by
$$a\vert \psi_{k,0}>=0\qquad {\overline a}\vert \psi_{k,k}>=0.\tag 23b$$
For $\nu_0$ assuming the above value, $\lambda_n$ now reads as
$$\lambda_n =(-1)^l \tan (\epsilon /2)[n]_{q^{1/2}}[k+1-n]_{q^{1/2}}.\tag 24 $$
It is important to note that the eigenstates are, as yet, not normalised.
To do so, we first compute the norm:
$$\align<\psi_{k,n}\vert\psi_{k,n}> &=<\psi_{k,n}\vert
{\overline a}\vert\psi_{k,n-1}>\\
&={\overline {\lambda}}_n<\psi_{k,n-1}\vert \psi_{k,n-1}> \\
&=\left(\prod\limits_{m=1}^n
\lambda_m \right)<\psi_{k,0}\vert\psi_{k,0}> \tag 25
\endalign$$
where we have replaced ${\overline {\lambda}}_m$ by $\lambda_m$
since it is real.
Now to ensure that this be positive definite, we set
$<\psi_{k,0}\vert \psi_{k,0}>=1$ and require that each
$\lambda_m$ term in the product be positive. Here the sign
is dictated by the term $(-1)^l\tan (\epsilon /2)$ (see (24))
and can be made positive if we
choose $l$ to be such that $(-1)^l\tan (\epsilon /2)>0$. For
instance, if $0<\epsilon < \pi /2$ then $l$ can be chosen to be
even.
\par
With this we can write the normalized states as
$$\vert k,n>={{\vert\psi_{k,n}>}\over{<\psi_{k,n}\vert \psi_{k,n}>^{1/2}}}
\qquad 0\le n\le k \tag 26$$
with
$$\align
&{\overline a}\vert k,n> = ((-1)^l\tan (\epsilon /2))^{1/2}
([n+1]_{q^{1/2}}[k-n]_{q^{1/2}})^{1/2}\vert k,n+1>\tag 27a\\
&a\vert k,n> = ((-1)^l\tan (\epsilon /2))^{1/2}
([n]_{q^{1/2}}[k+1-n]_{q^{1/2}})^{1/2}\vert k,n-1>\tag 27b\\
&N \vert k,n> =\left( \frac {-k-1}{2} +\frac
{(2l+1)\pi}{2\epsilon}+n\right) \vert k,n>.\tag 27c \endalign $$
Here it should be noted that the
representations are valid only for $\epsilon \ne 0$ or equivalently
$q\ne 1$. This is important, as one would otherwise be led to an
undeformed oscillator algebra with
finite dimensional representations. This clearly would have been a
paradox as it is well known that there are no such
representations for the usual oscillator algebra.
\par
On further examination, an even more interesting fact surfaces.
To see this, we first relabel the parameters $k$ and $n$ by
$k=2j$ and $n=j+m$. Then by
replacing $\vert k,n>\to \vert j,m>$, we have
$$\align
&((-1)^l\cot (\epsilon /2))^{1/2}{\overline a}\vert j,m>=
([j-m]_{q^{1/2}}[j+m+1]_{q^{1/2}})^{1/2}\vert j,m+1>\tag 28a\\
&((-1)^l\cot (\epsilon /2))^{1/2}a\vert j,m>=
([j+m]_{q^{1/2}}[j-m+1]_{q^{1/2}})^{1/2}\vert j,m-1>\tag 28b\\
&\left( N+\frac {1}{2}-\frac {(2l+1)\pi}{2\epsilon}\right) \vert j,m>=m\vert
j,m>\tag 28c
\endalign$$
where some terms have been transposed to the left hand side and $m$
and $j$ now assume
the values $-j\le m \le j$, $j=0,1/2,1,3/2,...$.
What is surprising here is that
the vectors on the right hand side
are precisely those of ${\text {su}}_{q^{1/2}}(2)$.
This immediately suggests that
$$\align
&J_+=((-1)^l\cot (\epsilon /2))^{1/2}~{\overline a}\tag 29a\\
&J_-=((-1)^l\cot (\epsilon /2))^{1/2}~{a}\tag 29b\\
&J_0=N +{{1}\over{2}}-{{(2l+1)\pi}\over{2\epsilon}} \tag 29c
\endalign$$
which essentially
establishes a linear relationship
between the generators of ${\Cal H}_q(1)$ and those of
${\text {su}}_{q^{1/2}}(2)$. It should be noted that
the relationship not only breaks down at $\epsilon =0$, but also at
the $\epsilon =p\pi$ and $\epsilon =(2p+1)\pi/2$ where $p\in {\bold
Z}$.
Now the identification has been made through the use of a
representation. One might be prompted to ask  whether they hold as
operator identities. To see if they do, all one needs to show is
that the commutation relations of one algebra reproduces the other
under (29). One can easily verify that this is indeed the case.
In addition the Hopf structure of ${\Cal H}_q(1)$ also induces through
the above relationship,
the Hopf structure of ${\text {su}}_{q^{1/2}}(2)$. This
implies that for $\vert q\vert =1$, modulo values for which (29)
is ill-defined, the Hopf $*$-algebra of ${\Cal H}_q(1)$ with
non-standard conjugation is equivalent to ${\text {su}}_{q^{1/2}}(2)$.
\par
Next,
let us consider the representations of ${\Cal H}_q(1)$ under the
involution provided by (13).
\plainfootnote {$^{\dag}$} {In the following we choose the
involution with the negative sign in (13).}
As in the previous case we
start with representations (18) with $\lambda_0=0$. Then under the
standard prescription
$$<\psi_n\vert N^+\vert \psi_n>={\overline {<\psi_n\vert N \vert
\psi_n>}}\tag 30$$
one is led to the condition
$$\nu_0-{\overline {\nu_0}}=i\frac {(2l+1)\pi}{\epsilon}\tag 31$$
which implies that
$$\nu_0=\mu_0+i\frac {(2l+1)\pi}{2\epsilon}\qquad \mu_0\in {\bold R}.
\tag 32$$
Let us see whether it is possible to truncate the states. The
existence of $\vert \psi_k>$ with ${\overline a}\vert \psi_k>=0$
now implies that
$$\cosh (\epsilon (\mu_0+k+1))=\cosh (\epsilon \mu_0)\tag 33$$
which is obtained in a manner analogous to (20a).
This in turn suggests that
$$\mu_0+k+1=\pm u_0\tag 34$$
which provides us with a possible solution
$k=-2\mu_0-1$. Since $\mu_0$ is arbitrary we can assign to it, the
values $-1/2,-1,-3/2,...$ which sets $k$ to $0,1,2,...$ as
required. Thus if we have
$$\nu_0= \frac{-k-1}{2}+i\frac{(2l+1)\pi}{2\epsilon}\qquad l\in
{\bold Z}\tag 35$$
then a truncation is certainly possible.
Under this condition $\lambda_n$ can be expressed as
$$\lambda_n=i(-1)^{l+1}\tanh (\epsilon /2)[n]_{q^{1/2}}
[k+1-n]_{q^{1/2}}\tag 36$$
with $0\le n \le k$ and $k=0,1,2,...$ The norms can be computed as in
(25) and are given by
$$<\psi_n \vert\psi_n>= \prod_{m=1}^{n}
\left( (-1)^{l+1}\tanh (\epsilon/2) [m]_{q^{1/2}}[k+1-m]_{q^{1/2}}\right).\tag
37 $$
To ensure positive definitness we choose $l$ to be even for
$\epsilon <0$ and odd for $\epsilon >0$.
Then using the form of (26) for the normalized states  we have
$$\align
&{\overline a}\vert k,n> = ((-1)^{l+1}\tanh (\epsilon /2))^{1/2}
([n+1]_{q^{1/2}}[k-n]_{q^{1/2}})^{1/2}\vert k,n+1>\tag 38a\\
&a\vert k,n> =i((-1)^{l+1}\tanh (\epsilon /2))^{1/2}
([n]_{q^{1/2}}[k+1-n]_q)^{1/2}\vert k,n-1>\tag 38b\\
&N \vert k,n> =\left( \frac {-k-1}{2} +i\frac
{(2l+1)\pi}{2\epsilon}+n\right)\vert k,n>.
\tag 38c
\endalign $$
It should be noted that
the above representations also break down at $\epsilon =0$.
Now these are very similar to (27) with the exception
of some complex factors in the terms. In view of this similarity, we are
again prompted to make the substitution
$k=2j$ and $n=j+m$ which yields
$$\align
&((-1)^{l+1}\coth (\epsilon /2))^{1/2}{\overline a}\vert j,m>=
([j-m]_{q^{1/2}}[j+m+1]_{q^{1/2}})^{1/2}\vert j,m+1>\tag 39a\\
&-i((-1)^{l+1}\coth (\epsilon /2))^{1/2}a\vert j,m>=
([j+m]_{q^{1/2}}[j-m+1]_{q^{1/2}})^{1/2}\vert j,m-1>\tag 39b\\
&\left( N+\frac {1}{2}-i\frac {(2l+1)\pi}{2\epsilon}\right) \vert j,m>=m\vert
j,m>\tag 39c
\endalign$$
with the implication that
$$\align
&J_+=((-1)^{l+1}\coth (\epsilon /2))^{1/2}~{\overline a}\tag 40a\\
&J_-=-i((-1)^{l+1}\cot (\epsilon /2))^{1/2}~{a}\tag 40b\\
&J_0=N +{{1}\over{2}}-i{{(2l+1)\pi}\over{2\epsilon}}. \tag 40c
\endalign $$
After affirming that the relationship holds at the
algebraic level, we again conclude that ${\Cal H}_q(1)$ together
with the conjugation (13) is
equivalent to ${\text {su}}_{q^{1/2}}(2)$ for $\epsilon
\ne 0$ but now with the standard conjugation.
\eject
\noindent
{\bf REFERENCES}
\ref 1 Faddeev, L.D., {\it Les Houches Lectures} 1982
(Elsevier,Amsterdam, 1984).
\spa {~} Kulish, P.P. and Sklyanin, E.K., {\it Lecture Notes in Physics} Vol.
151 Springer, Berlin, 1982.
\ref 2 Yang, C.N., {\it Phys. Rev. Lett.}{\bf 19}, 1312 (1967);
\spa {~} Baxter, R.J., {\it Exactly Solved Models in Statistical
Mechanics}, New York: Academic, 1982.
\ref 3 Zamolodchikov, A. and Zamolodchikov, Ab., {\it Ann. Physics}
{\bf 120}, 252 (1979);
\spa {~} de Vega, H., {\it Int. J. Mod. Phys.}{\bf 4}, 2371 (1989).
\ref 4 Moore, G. and Seiberg, N., {\it Nucl. Phys.}B{\bf 313},
16 (1989); {\it Commun. Math. Phys.}{\bf 123}, 177 (1989).
\ref 5 Drinfeld, V.G., {\it Proc. Intern. Congress of
Mathematicians}, Berkley Vol. 1, 1986, p 798.
\ref 6 Macfarlane, A.J., {\it J. Phys.}{\bf A22}, 4581 (1989).
\ref 7 Biedenharn, L.C., {\it J. Phys.} {\bf A22}, L873 (1989).
\ref 8 Yan, H., {\it J. Phys.}{\bf A23}, L1155 (1990).
\ref {9} Floreanini, R. and Vinet, L., {\it Lett. Math. Phys.}{\bf
22}, 45 (1991).
\ref {10} Masuda, T., Mimachi, K., Nakagami, Y., Noumi, M., Saburi,
Y. and Ueno, K., {\it Lett. Math. Phys.}{\bf 19}, 187 (1990).
\ref {11} Lukierski, J., Nowicki, A. and Reugg, H., {\it Phys. Lett. B}
{\bf 271}, 321 (1991).
\ref {12} Fr{\"o}hlich, J. and Kerler, T., {\it Lecture Notes in
Mathematics} Vol 1542, Springer-Verlag, Berlin, 1993
\bye